\documentclass[a4paper,11pt]{article}
\usepackage{amsmath}
\usepackage{amsfonts}
\usepackage{amscd}
\usepackage{amssymb}

\newcommand{\al}{{\alpha}}

\newcommand{\om}{{\omega}}

\newcommand{\ra}{\rightarrow }

\newcommand{\MS}{{\medskip}}
\newcommand{\SK}{{\smallskip}}
\newcommand{\BS}{{\bigskip}}
\newcommand{\NI}{{\noindent}}

\newcommand{\QED}{\hfill {\bf QED}\medskip}

\newtheorem{theorem}{Theorem}[section]
\newtheorem{thm}[theorem]{Theorem}
\newtheorem{cor}[theorem]{Corollary}

\newtheorem{lemma}[theorem]{Lemma}
\newtheorem{prop}[theorem]{Proposition}
\newtheorem{ex}[theorem]{Example}

\numberwithin{equation}{section}
\newtheorem{nta}[theorem]{Notation}

\title{\bf Symplectic forms invariant under free circle actions on
4--manifolds}
\author{Bogus\l\/aw Hajduk $^*$ and Rafa\l\ Walczak \thanks{Both authors
were partially supported by Grant 2 P03A 036 24 of Polish
Committee of Sci. Research} }

\date{\today}

\begin{document}

\maketitle

\begin{abstract}
\NI Let $M^4$ be a smooth compact manifold with a free circle
action generated by a vector field $X.$ Then for any invariant
symplectic form $\om$ on $M$ the contracted form $\iota_X\om$ is
nondegenerate. Using the map $\om \mapsto \iota_X\om$ and the
related map to $H^1(M \slash S^1,\mathbb R)$ we study the topology
of the space $S_{inv}$ of invariant symplectic forms. In
particular we give a description of $\pi_0S_{inv}(M^4)$ in terms
of the unit Thurston ball.

\smallskip


\noindent {\bf AMS classification(2000)}: 53D05, 57S25

\end{abstract}

\bigskip

\bigskip

\section{Introduction}

Consider a smooth compact manifold $M.$ A symplectic form on $M$
is a closed nondegenerate differential $(C^\infty )$ 2--form
$\omega.$ A manifold equipped with a symplectic form is called
symplectic manifold or we say that the form yields a symplectic
structure on the manifold. We will always assume $M$ oriented and
symplectic forms compatible with the orientation, which means
$\omega^n > 0.$ The space of all symplectic forms on $M$ with
$C^\infty$ topology is, if nonempty, an infinitely dimensional
space. The main motivation for this paper was to understand the
homotopy of this space. The simplest question here is whether it
is connected or what is the set of connected components.  We say
that two symplectic forms are homotopic if they are in the same
component of the space of symplectic forms and that they are
isotopic if they are in the same component of the space of
symplectic forms in a given cohomology class in $H^2(M,\mathbb
R).$ Thus two symplectic forms $\om_0, \om_1$  are homotopic if
they can be joined by  a continuous  path of symplectic forms
$\om_t, t\in[0,1]$ and isotopic if there is  a path such that the
cohomology class of $\om_t$ is constant. Thus we can ask under
what conditions two given symplectic forms are homotopic
(isotopic).

The only easy case is dimension 2, when a symplectic form is an
area form. Two area forms are homotopic if and only if they give
the same orientation and any two cohomologous symplectic forms are
isotopic. In higher dimensions there  are very few cases for which
some answers  are known. In dimension 4  the problem is open even
for  the torus $T^4.$ A useful and often calculable invariant is
provided by almost complex structures, i.e.  complex structures on
the tangent bundle $TM.$ Thus an almost complex structure is a
bundle automorphism $J$ of $TM$ such that $J^2=-Id.$  For any
symplectic form $\omega$ there exist  almost complex structures
$J$ on $M$ tamed by $\omega $ in the sense that $\omega (V,JV)>0$
for any nonzero  tangent vector $V.$ The space of almost complex
structures tamed by $\omega$ is infinite dimensional but
contractible and for any continuous path $\omega_t$ there exists a
path $J_t$ such that $J_t$ is tamed by $\omega_t.$ Thus Chern
classes of $(TM, J)$ are equal for $J's$ tamed by  homotopic
symplectic forms. Again, some simple questions  are still open.
For example for $T^4$ there are homotopy class of almost complex
structures for which it is not known whether they are tamed by a
symplectic form.

When one studies the case of the torus, this is  easy to find that
the space of harmonic (with respect to the canonical flat
Riemannian metric) symplectic forms compatible with fixed
orientation is connected. Such forms have constant coefficients
and the nondegeneracy condition is a quadratic inequality which is
easy to deal with. However, there is a  larger class of symplectic
forms which is not so easy to understand and which can serve as
intermediate class between harmonic and general symplectic forms.
This is the space of those symplectic forms which are invariant
with respect to an action of the circle $S^1$ on the torus.

The primary aim of this paper is to study homotopy properties of
the space $S_{inv}(M)$ of invariant symplectic forms compatible
with orientation on a closed oriented 4--manifold $M$ with a free
action of the circle.

Some theorems about existence and structure of such forms and
manifolds were proved in \cite{Ba1,Bo,FMG}.

The basic observation is that a symplectic form $\om$ on $M$
invariant under a free action of a circle gives, under contraction
with the infinitesimal generator $X$ of the action, a
nondegenerate closed form $\iota_X \om .$ Thus there is also a
closed nondegenerate form having rational periods and this implies
that  $M$ fibers over the circle \cite{Ti}.
 Moreover, since
$\iota_X \om$ is  invariant and it vanishes on vertical bundle of
the fibration by circles $\pi : M\rightarrow M/S^1,$ there is a
closed nondegenerate form $\alpha$ on $M/S^1$ such that $\iota_X
\om = \pi^*\alpha .$ Thus  the orbit space $M/S^1$ also fibers
over $S^1.$ This is a severe topological restriction if $M$ is
closed. In dimension four it implies that $M$ is either aspherical
or diffeomorphic to $S^2 \times T^2.$

There is also a restriction for the cohomology class of $\alpha :$
we always have  $c_1(\pi )\cup [\al]=0, $ where $c_1(\pi )$ is the
first Chern class of the circle fibration $\pi$ (see lemma
\ref{bu}). The last condition was proved in \cite{Bo}, by a kind
of inflation trick, to be sufficient to existence of an invariant
symplectic form determining a given nondegenerate and closed form
$\alpha$ (see \ref{sim}).

In many instances the questions about the topology of the space of
symplectic forms stack over the problem of weird properties of
some natural maps. For example, we have the map to $ H^2(M,\mathbb
R)$ which associate to each symplectic form its cohomology class
or the map to moduli space of symplectic forms. If this map had
decent properties, as for example the lifting path property, then
one can apply the standard machinery of algebraic topology to
study topology of symplectic forms. Unfortunately, examples as in
\cite{MD1} show that this is not true in general, at least in
dimensions greater than 4.

In the invariant  case the map $\omega \mapsto \al$ provides
dimensional and degree reduction so the problems on invariant
symplectic forms can be expressed in terms of 1--forms (cf.
Introduction in \cite{MT}). We show that in dimension four this
map as well as the induced map to $H^1(M,\mathbb R)$ have quite
nice properties. Our main result is that it is a homotopy
equivalence of the space of invariant symplectic forms
$S_{inv}(M)$ compatible with the orientation to a subspace of
nondegenerate closed 1--forms (Theorem \ref{im}). There are three
principal ingredients of the proofs. First is the Blank --
Laudenbach theorem which says that any two closed cohomologous
nondegenerate 1--forms on a 3--manifold are isotopic. Secondly,
Thurston norm  is used to decide which cohomology classes in $
H^1(M,\mathbb R)$ are represented by nondegenerate closed forms.
Finally, to pass from 1--form to symplectic forms we use the
inflation formula \ref{sim}. The main technical step is that any
path of 1--forms satisfying a necessary condition (\ref{bu}) lifts
to a path of symplectic forms (Lemma \ref{lift}). As a corollary
we get that any two symplectic forms on $T^4$ invariant with
respect to any free circle action and determining the same
orientation are homotopic through invariant symplectic forms.

Next we use this result to give sample calculations of $S_{inv}.$
For the torus $T^4$ we prove that $S_{inv}$ is simply connected
and for every nontrivial circle bundle over $T^3$ we have
$\pi_1(S_{inv}) = \mathbb Z.$ We show also an example of
3--manifold and a cohomology class of degree one such that the
space of nondegenerate forms representing this class has the
homotopy type of the circle (i.e. the Blank -- Laudenbach theorem
do not extend to higher homotopy groups). Another interesting
space which appears in our considerations is the space of
automorphisms of a codimension 1 foliation on $N^3.$ We use
calculations of that space for the foliation by (compact) fibers
of a fibration over $S^1,$ but possibly computations for $S_{inv}$
can give some information on automorphisms of foliations with
noncompact leaves.

In the last section, we study from a similar point of view  the
correspondence between symplectic forms and almost complex
structures in the invariant case. It is known that for symplectic
4--manifold there exist homotopy classes of almost complex
structures not tamed by any symplectic form. Nontrivial examples
of that kind were given in \cite{CLO} be twisting a given $J$ in a
disc and calculating that the {\it Seiberg --Witten} invariant
changes, which contradicts the Taubes theorem  \cite{TA1,TA2}.
Such twist can be done along an orbit of any linear circle action
in invariant way, but the Seiberg -- Witten invariant does not
change. It is not known whether resulting almost complex
structures are tamed by a (maybe not invariant) symplectic form.
Nevertheless, we give a sufficient condition for an invariant $J$
to be homotopic to one which is tamed by an invariant symplectic
form. We prove also the following lifting property. Suppose
$J_0,J_1$ are tamed by symplectic structures $\om_1,\om_2$ and
there is a continuous path $J_t, t\in [0,1]$ of invariant almost
complex structures. In general it is not possible to lift the path
to such a path $\om_t$ that $\om_t$ tames $J_t,$ but we show that
there exists another path with the same ends which lifts. We give
also a sufficient condition for an invariant $J$ to be homotopic
to one which is tamed by an invariant symplectic form.

Our results can be possibly continued in two directions. One is
the natural extension to locally free actions. In that case the
orbit space is an orbifold and the inflation trick does work. The
problem is that Blank -- Laudenbach theorem and the Thurston norm
theory are known only for smooth manifolds. Another potentially
interesting extension is to the space of forms that are invariant
under the flow of a nowhere vanishing vector field. Then the
topological restrictions are  relaxed ($M$ still fibers over the
circle, but need not be aspherical). When $M$ has no locally free
circle action, it can be thought as a substitute of invariant
forms. Even if we have a free circle action on $M,$ such space may
be interesting as an intermediate space between invariant and
general symplectic forms.

\section{Homotopy classification of symplectic forms invariant
under free circle actions on 4 -- manifolds}\label{sym}

Throughout the article $M$ will denote a compact oriented smooth
manifold with e smooth free action of $S^1$ and $N \cong M \slash
S^1.$ The space of invariant symplectic forms consistent with the
given orientation will be denoted by $S_{inv}.$

Let $\pi:M \ra N$ denote the principal $S^1$ -- fibration given by
the action. By $\al$ we will denote the non\-dege\-ne\-rate
(=no\-where va\-ni\-shing) closed 1--form $\alpha$ on $N$
satisfying
$$\pi^*\alpha = \iota_X\omega, $$
where $X$ is the infinitesimal generator of the action.

\begin{lemma}\label{bu} {\bf \cite{Bo}} If $\omega \in S_{inv}$, then
$$[\alpha] \cup c_1(\pi) = 0.$$
\end{lemma}

{\bf Proof}. Take any connection form $\eta \in \Omega^1(M,\mathbb
R)$. By Chern -- Weyl there exists a closed 2 -- form $c_1 \in
H^2(N,\mathbb R)$ such that $\pi^*c_1 = d\eta$ and $c_1$
represents $c_1(\pi) \in H^2(N^3,\mathbb R)$. There also exists a
unique 2 -- form $\beta' \in \Omega^2(N^3,\mathbb R)$ such that
$\omega - \eta \wedge \iota_X\omega = \pi^*\beta'$. By
differentiating both sides of the last equation we get $d\eta
\wedge \iota_X\omega = \pi^*d\beta'$ and this implies the lemma.
\QED

\MS

From now on we assume that $\dim M=4.$

\MS

Define the subspace $L \subset H^1(N^3,\mathbb R)$ by

\begin{eqnarray} L = \{ \alpha \in H^1(N^3,\mathbb R) \mid
\alpha \cup c_1(\pi) = 0\}. \end{eqnarray}

Let us also recall the notion of {\bf Thurston norm} \cite{Th}
(see also \cite{MT}) for a 3--manifold. If $N^3$ is a compact,
connected and oriented manifold without boundary then for any
compact oriented n--component surface $S = S_1 \sqcup \cdots
\sqcup S_n$ embedded in $N$ define

\begin{eqnarray} \chi_-(S) = \sum_{\chi(S_i) <0} \arrowvert
\chi(S_i) \arrowvert. \end{eqnarray}

The {\bf Thurston norm} on $H_2(N^3,\mathbb Z)$ and, by the
Poincare duality, on $H^1(N^3,\mathbb Z)$ is given by

\begin{eqnarray} {\Arrowvert \phi \Arrowvert}_T = \inf \{ \chi_-(S)
\mid [S] = \phi  \}.\end{eqnarray}

The Thurston norm can be extended linearly to $H^1(N^3,\mathbb
R)$. Let $B_T = \{\phi \ : \  {\Arrowvert \phi \Arrowvert}_T \leq
1 \}$ denote the unit ball in the Thurston norm. It is a (possibly
noncompact) polyhedron in $H^1(N^3,\mathbb R)$. Suppose $\phi' \in
H^1(N^3,\mathbb Z)$ is represented by a fibration $N^3 \rightarrow
S^1$. Then $\phi'$ is contained in the open cone $\mathbb R_+
\cdot F$ over a top -- dimen\-sional face $F$ of the Thurston norm
ball $B_T$. In this case we say $F$ is a {\it fibered face} of the
Thurston norm ball. The Thurston norm can be also defined if
boundary of $N^3$ is a union of tori.

Assume now that $N^3$ fibers over the circle, $b_1=m$ and it is
not $S^2 \times S^1$. In this case we have that $L \cap \mathbb
R_+ \cdot F$ has a homotopy type of a point, or a sphere $S^{m-1}$
if $c_1=0$ in $H^2(N,\mathbb R),$ and $S^{m-2}$ when $c_1 \neq 0$
provided that Thurston norm vanishes identically. If $N \cong S^2
\times S^1$, then total space of any circle fibration over $N$
happens to be symplectic if and only if $c_1=0$ and is
diffeomorphic $S^2 \times T^2$; furthermore we have that $L \cap
\mathbb R_+ \cdot F \cong \mathbb R \backslash \{0\}$ has a
homotopy type of a two--point set.

\MS

\begin{nta} For a given 3--dimensional manifold $N$ let's
define by $\mathcal{N}$ the set of all closed and nondegenerate
1--forms, by $\mathcal{N}_L \subset \mathcal{N}$ these closed and
nondegenerate 1--forms whose cohomology classes are included in
$L$, by $\mathcal{N}_H \subset H^1(N,\mathbb R)$ the set of all
cohomology classes which are represented by a nondegenerate,
closed 1--form, and finally by $\mathcal{N}^\alpha \subset
\mathcal{N}$ the set of all nondegenerate, closed 1--form whose
cohomology class equals $\alpha \in H^1(N,\mathbb R)$.
\end{nta}

\MS

Later we will need the following.

\begin{lemma}\label{lab} Assume that closed three--dimensional manifold $N$
fibers over the circle and we are given a path $\{x_t\}_{t \in
[0,1]} \subset \mathcal{N}_H$. Then there is a path
$\{\alpha_t\}_{t \in [0,1]} \in \mathcal{N}$ such that
$[\alpha_t]=x_t$ for all $t \in [0,1]$. We can also assume that
$\alpha_0$ and $\alpha_1$ are prescribed.
\end{lemma}

{\bf Proof}. For each $t \in [0,1]$ there exists a path
$\alpha_t'$ with the desired property on some neighborhood
$(t-\epsilon,t+\epsilon)$ of $t$. Choose a finite subcover and
glue together these paths on their overlaps using the Blank --
Laudenbach \cite{BL}, which states that any two cohomologous,
nondegenerate and closed 1--forms on a three--dimensional closed
manifold are isotopic. Isotopy here can be understood like in the
case of symplectic forms in two ways. The first is where two
nondegenerate closed cohomologous 1--forms can be joined by a path
of nondegenerate closed cohomologous forms, while the second means
that we have a path of diffeomorphisms $\{\phi_t\}_{t \in [0,1]}$
such that $\phi_0=Id$ and $\phi_1^*\al_1=\al_0.$ These two notions
are equivalent due to Moser argument \cite{Mo}.

\QED

\begin{cor} The mapping $\phi_N: \mathcal{N} \rightarrow \mathcal{N}_H$
has the property of lifting any circle mappings.
\end{cor}
\smallskip

The basic problem that emerges now is to decide, in terms of
$\al_1,\al_2,$ when two symplectic invariant forms $\om_1,\om_2$
are homotopic in $S_{inv}.$  Recall that (see \cite{MT}) there
exist on a 3--dimensional manifold $N$ two nonhomotopic closed
1--forms $\al_1,\al_2$ and, on $N \times S^1$ two nonhomotopic
symplectic forms.

The main technical lemma of this section gives an answer to the
problem.

\begin{lemma}\label{lift}
If we are given a path $\{ \alpha_t \}_{t \in [0;1]} \subset
\mathcal{N}_L$, then there exists a path $\{ \omega_t \}_{t \in
[0;1]}$ in $S_{inv}$ such that $\pi^*\alpha_t = \iota_X\omega_t$
for every $t \in [0;1]$.
\end{lemma}

{\bf Proof}. it follows from a formula \cite{Bo} in the spirit of
the {\it inflation} trick of Thurston \cite{Th1} and McDuff
\cite{MD2}. It consists in enlarging the form along the foliation
determined by $\ker\alpha.$ For a given form $\alpha,$

\begin{eqnarray}\label{sim} \omega = \eta \wedge \pi^*\alpha + \pi^*(K \beta +
\phi)\end{eqnarray}

\noindent is an invariant symplectic form if $\eta \in
\Omega^1(M^4,\mathbb R)$ is a connection form, $\beta$ is a closed
2 -- form on $N^3$ such that $\alpha \wedge \beta$ is a volume
form on $N^3, \ d\phi = -c_1 \wedge \alpha$ and $K$ is
sufficiently large real number. Obviously $\omega$ satisfies
$\pi^*\alpha = \iota_X\omega$.

\smallskip

Existence of $\beta$ is well--known.

\begin{lemma}\label{bet} Let $N^n$ be a closed and oriented manifold. Assume
that a closed and nondegenerate 1--form $\alpha$ on $N$ is given.
Then there is a closed $(n-1)$--form $\beta$ such that $\alpha
\wedge \beta$ is a volume form on $M$. Equivalently, $\beta$ is
nondegenerate on leaves of the foliation defined by $\ker\alpha.$
\end{lemma}

We use Lemma \ref{bet} to construct a continuous path of closed
2--forms $\{\beta_t\}_{t \in [0,1]} \subset H^2(N,\mathbb R)$ such
that $\alpha_t \wedge \beta_t$ are volume forms for all $t \in
[0,1]$. Indeed, for each $t$ there is $\beta_t$ which is
appropriate for some neighborhood of $t$. As in the proof of lemma
\ref{lab} we glue together patches of $\beta_t$ on $U_t$ using the
simple fact that if $\alpha_t \wedge \beta_t$ and $\alpha_t \wedge
\beta_t'$ are volume forms, then $\alpha_t \wedge (s \beta_t
+(1-s)\beta_t')$ are volume forms for $s \in [0,1]$ as well.

Next choose a continuous path $\phi_t$ which satisfies the
condition $d\phi_t = -c_1 \wedge \alpha_t$. To finish notice that
for $K$ large enough
\begin{eqnarray}\omega_t = \eta \wedge \pi^*\alpha_t + \pi^*(K \beta_t +
\phi_t)\end{eqnarray} is an invariant symplectic form for all $t$.
\QED

In order to proceed we need the following simple lemma.

\begin{lemma}\label{con} If $\omega, \omega' \in S_{inv}$ and
$\alpha = \alpha'$, then $\omega, \omega'$ are homotopic in
$S_{inv}$. If we assume also that $[\om]=[\om'],$ then $\om$ is
isotopic to $\om'.$ The same holds in the parametric case, i.e.
for two continuous families of symplectic forms parameterized by a
topological space.
\end{lemma}

{\bf Proof}. Under our assumptions $\omega$ and $\omega'$ can be
joined in $S_{inv}$ by the convex combination. To see this write
$$\omega= \eta \wedge \pi^*\alpha + \pi^*\phi, \ \
\omega'= \eta \wedge \pi^*\alpha + \pi^*\phi'$$ for some 2--forms
$\phi,\phi'$ on $N$. Then $$(\lambda \omega+(1-\lambda)
\omega')^2=(\eta \wedge \pi^*\alpha +\pi^*(\lambda
\phi+(1-\lambda) \phi'))^2=$$ $$=2(\lambda \omega^2+(1-\lambda)
\omega'^2)>0.$$

\QED

\MS

\NI More generally for cohomological $\om,\om'$ we have the
following.

\begin{lemma}\label{in} Any two invariant,
cohomological and symplectic forms $\omega$ and $\omega'$ are
isotopic.
\end{lemma}

{\bf Proof}. It follows from Blank--Laudenbach \cite{BL} that
there exists an isotopy $\{\phi_t\}$ between closed 1 -- forms
$\alpha,\alpha'$ on $N^3$ obtained from $\omega,\omega'$ such that
$\phi_1^*\al'=\al$ and $\phi_0=Id.$ The path $\phi_t$ lifts to a
path $\widetilde{\phi_t}$ on $M$ such that
$(\widetilde{\phi_t})_*X=X.$ Then
$$\iota_X(\widetilde{\phi_1})^*\om'=
(\widetilde{\phi_1})^*\iota_X\om'=(\widetilde{\phi_1})^*\pi^*\al'=
\pi^*\phi_1^*\al'=\pi^*\al=\iota_X\om$$ and
$\om,(\widetilde{\phi_1})^*\om'$ are again in the same cohomology
class. To finish the proof apply \ref{con}.  \QED

\SK

Now we are in position to prove our main results.

\begin{thm} Let $M^4$ be a closed manifold equipped with a free
action of the circle $S^1$. Assume that we are given a path $\{x_t
\}_{t \in [0,1]} \subset H^2(M^4,R)$ such that for each $t \in
[0,1]$ there exists an invariant symplectic form $\omega_t'$ in
the class that $x_t.$ Then there exists a path $\{\omega_t \}_{t
\in [0,1]}$ of symplectic and invariant forms such that
$[\omega_t]\in x_t$. Moreover, the forms $\omega_0,\omega_1$ can
be prescribed.
\end{thm}

{\bf Proof}. For every $t \in [0,1]$ there exists $\epsilon >0$
such that in the closed interval $[t-\epsilon,t+\epsilon]$ we have
a path of symplectic and invariant forms satisfying the thesis of
our theorem.

\noindent By the standard compactness argument there exists an
increasing sequence $0 = t_0 < t_1 < \ldots < t_n =1$ such that in
the closed interval $[t_i,t_{i+1}]$ the given path lifts.
Furthermore forms standing over $t_i$ can be glued with the help
of the isotopy given by Lemma \ref{in}. \QED

\BS

\noindent In higher dimensions our arguments do not work. Let us
recall the example given by McDuff.

\noindent On $M = S^2 \times S^2 \times T^2$ we have a family of
symplectic forms \begin{eqnarray}\tau_\lambda=\lambda \sigma_0
\times \sigma_1 \times \sigma_2, \ \ \ \lambda \geq
1,\end{eqnarray} \noindent where $\sigma_i$ are area forms of
total area equal to 1. If we take diffeomorphism $\phi : M
\rightarrow M$ given by the formula
\begin{eqnarray}\phi(z,w,s,t)=(z,\rho_{z,t}(w),s,t)\end{eqnarray} \noindent where
$\rho_{z,t}:S^2 \rightarrow S^2$ is the rotation around the axis
through $z$ through the angle $2\pi t$, then it was shown
(\cite{MD1}) that $\tau_\lambda$ and $\phi^*\tau_\lambda$ are
isotopic when $\lambda >1$ but are not isotopic when $\lambda =
1$.

\noindent Let us define on $M$ a free circle action by the formula
\begin{eqnarray}\theta(z,w,s,t)=(z,w,s \theta,t), \ \ \ \theta \in S^1.\end{eqnarray}
\noindent Then $\phi$ commutes with the action, and $\tau_\lambda$
is invariant when $\sigma_2$ is chosen to be invariant.

\SK

\noindent Observe that we have \begin{eqnarray}
\iota_X\tau_\lambda = \iota_X\phi^*\tau_\lambda =dt,
\end{eqnarray} if $X=\frac{\partial}{\partial s}$ denotes
the infinitesimal generator of the action.

\NI We can easily check (see \cite{MD1}) that $\tau_\lambda$ and
$\phi^*\tau_\lambda$ are isotopic through invariant symplectic
forms for $\lambda>1$, and that $\tau_1,$ $\phi^*\tau_1$ are not
isotopic.

\NI It follows that Lemma \ref{con} does no longer work in
dimensions higher than 4.

\hfill $\Box{}$

\MS

The mapping $\omega \mapsto \alpha$ satisfies a stronger property
then a path lifting property. Namely, this map is a homotopy
equivalence.

\begin{thm}\label{im} The map \begin{eqnarray}\Phi:S_{inv}\rightarrow
\mathcal{N}_L,\end{eqnarray} defined by $\omega \mapsto \alpha$,
is a homotopy equivalence.
\end{thm}

{\bf Proof}. We will show that the map $\Phi_*: \pi_nS_{inv}
\rightarrow \pi_n\mathcal{N}_L$ is an  isomorphism for every $n$.
Since both spaces have the homotopy type of CW--complexes, we get
the conclusion by the Whitehead theorem.

Let us first consider the case where $n=0$. For $n
> 0$ the reasoning is a direct extension of the proceeding
case. We need only to show that given $\om_0,\om_1$ a homotopy
between $\al_0,\al_1$ rises to a homotopy between $\om_0$ and
$\om_1.$ The idea of the proof is to take any lift $\om_t$ given
by (\ref{lift}), and some partial lifts of $\al_t$ over some
neighborhoods of $0$ and $1$ which extend given forms $\om_0$ and
$ \om_1$ respectively. Using the convex combination we can patch
together these paths by the (parametric version) of \ref{con} to
get one such that at the ends we get the prescribed forms. We
leave the details to the reader.

\smallskip

For the general case, to see that $\Phi_*$ is onto assume we have
a map $\phi_1:(S^n,\star) \rightarrow (\mathcal{N}_L,\star)$. We
have to prove the existence of lifting of $\phi_1$ to $S_{inv}$.
There exists a map $\phi_2:S^n \rightarrow \text{cl}(\Omega^2N)$
such that $\phi_1(\theta) \wedge \phi_2(\theta)$ is a volume form
on $N^3$ for every $\theta \in S^n$. To construct $\phi_2$ follow
the construction in lemma \ref{lift}. First find a finite open
cover of $S^n$ such that appropriate condition is satisfied on
each element of the covering, then glue $\beta_\theta$ using a
partition of unity subordinated to the covering. To finish the
proof use the (parametric version) of the Lemma \ref{con}.

Proof that $\Phi_*$ is an injection is very similar. \QED

\smallskip

\begin{cor}\label{tor} Let $\om$ be a symplectic form on $T^4$
invariant under a linear action of $S^1.$ Then $\omega$ is
homotopic to the standard structure $\omega_{st}=dx_1 \wedge
dy_1+dx_2 \wedge dy_2$. Furthermore, $\om$ is isotopic to some
constant coefficient form.
\end{cor}

{\bf Proof}. If $S^1$ acts linearly on $T^4$, then the quotient
space is diffeomorphic to $T^3$. The Thurston norm equals $0$
everywhere, so we have $\mathcal{N}_H \cong H^1(T^3,\mathbb R)
\backslash \{0\}$. It yields that $\mathcal{N}_H = \mathcal{N}$ is
connected, so due to theorem \ref{im} we prove the corollary. \QED

\NI {\bf Remark}. By a linear action we mean the action equivalent
to an action of a subgroup of the torus considered as the Lie
group. It follows from the Poincare conjecture then any free
action of the circle is linear.

\section{Symplectic forms and diffeomorphisms}\label{Sydi}

In the sequel we shall use the fact that the evaluation maps are
fibrations in some cases. Firstly, if we have a topological group
acting transitively on a space and there are locals sections of
the {\it ev}, then it is locally product. Secondly, the evaluation
map of $\mathcal{M}(X,Y) \ra Y \ : \ f \mapsto f(x_0),$ given by
the base point $x_0 \in X$ is Hurewicz if $\mathcal M$ is a
component of all continuous maps from $X$ to $Y$ and $X,Y$ are
both compact CW -- complexes. For details consult \cite{Sp}.

 The homotopy properties of $S_{inv}$ can be calculated,
 in simple cases, using the action of $Diff_0$ (= diffeomorphisms
 isotopic to the identity) on forms. Consider the evaluation map
\begin{eqnarray}ev_\alpha : Diff_0N
 \longrightarrow \mathcal{N}^\alpha,\end{eqnarray}
 defined by
\begin{eqnarray}ev_\alpha(\phi)=\phi^*\alpha.\end{eqnarray}
 The map $ev_\alpha$ is onto due to \cite{BL}, so
\begin{eqnarray}\label{fi}F {\buildrel \kappa \over
\longrightarrow} Diff_0N {\buildrel ev_\alpha \over
\longrightarrow} \mathcal{N}^\alpha,\end{eqnarray} is a locally
trivial fibration with fiber  $F$ equal to the space of those
diffeomorphisms of $N$ which are isotopic to the identity and
preserve the form $\alpha$ or, equivalently, the foliation defined
by $\alpha$. Thus $\mathcal{N}^\alpha$ is contractible if and only
if the inclusion $\kappa : F \rightarrow Diff_0N$ is a homotopy
equivalence.

We will give calculation for the torus $T^4$ and $S^1$--bundle
over $T^3.$ These cases are particularly simple because the
inclusion $T^k \hookrightarrow Diff_0T^k,$ given by translations,
is a homotopy equivalence for $k \leq 3.$

\MS

\noindent For the 3--dimensional torus we claim that
$\mathcal{N}^\alpha$ is contractible  for any closed and
nondegenerate 1--form $\alpha$  defining a fibration, i.e., a form
of the form $\tau^*dt$ for a fibration $\tau$ and the standard
1--form $dt$ on $S^1.$ This property depends only on the
cohomology class of $\alpha,$ since it is equivalent to that the
cohomology class of $\alpha$ belongs to the positive cone over
integer classes. We will call any cohomology class {\bf fibrated}
if it is represented by such a form.

Any fibration of $T^3$ over the circle is product, so we can
assume $\alpha=dx_1$. Then $F$ consists of  fiberwise
diffeomorphisms which covers rotations on $S^1$,  uniquely
described by two parameters: an element of the space $Diff_0T^2
\times \Omega_0Diff_0T^2$ and the rotation angle $\theta \in S^1$,
where $\Omega_0$ denotes the space of nullhomotopic loops. We have
the following sequence of homotopy equivalences:

$$Diff_0T^2 \times \Omega_0Diff_0T^2\times S^1 \cong
T^2 \times \Omega_0T^2\times S^1 \cong T^2 \times \{\star\}\times
S^1 \cong T^3.$$

This follows from $Diff_0T^2 \cong T^2$ and the asphericity of
$T^2.$

In the diagram

$$
\begin{CD}\label{ull}
             T^3           @>{\hbox{Id}}>>   T^3 \\
             @V{}VV              @V{ }VV   \\
             F           @>{\kappa}>>     Diff_0T^3 \\
\end{CD}
$$

\noindent  vertical maps induced by inclusions, are homotopy
equivalences, and so is $\kappa$. Hence $\mathcal{N}^\alpha$ is
contractible for a fibrated class.

\MS \noindent We can now derive the following corollary of
\ref{im}:

\begin{cor}\label{mwk} For any linear free action of $S^1$ on
$T^4,$ $S_{inv}$ is simply connected.
\end{cor}

{\bf Proof}. Since in our case $\mathcal{N}_L=\mathcal{N}$, so it
is enough to prove that $\pi_1\mathcal{N}=0.$ We will deform
continuously a given loop $\{\alpha (t),\ t \in S^1\}$ of closed
and nondegenerate 1--forms to the space of forms with constant
coefficients.   For every $t_0 \in S^1,$ $\alpha(t_0)$ is isotopic
to a form of constant coefficients. Moreover, there exists a
neighborhood $U$ of $t_0$ and an extension of the isotopy to a
mapping
$$\alpha (s,t) : U \times [0,1] \longrightarrow  \mathcal{N}$$ such
that $\alpha_{s,t}$ is an isotopy for fixed $t\in U$ and
$\alpha(1,t)$ is a path of constant coefficients forms. This is
given by $\alpha(s,t_0)+\delta (t)-\delta(t_0),$ where $\delta(t)$
is the form of constant coefficient in the cohomology class of
$\alpha (t).$ If $U$ is small enough, then the convex combination
makes a isotopy between $\alpha (t)$ and $\alpha(t_0)+\delta
(t)-\delta(t_0).$

Considering a finite cover by such neighborhoods we come to a
finite partition of $S^1$ by $t_1,\ldots,t_k \in S^1$ such that on
each interval $[t_i,t_{i+1}]$
 we have the required homotopy.
Over each point $t_i$ there are two isotopies joining
$\alpha(t_i)$ with the form of constant coefficients, thus we see
a loop in $\mathcal{N}^{\alpha(t_i)}.$ If we perturb the loop by a
form $u$ of constant coefficients, then we get a loop in
$\mathcal{N}^{\alpha +u}.$ But we can choose $u$ such that
$(\alpha +u)$ is fibrated class and thus $\mathcal{N}^{\alpha +u}$
is contractible. If $u$ is small enough, the perturbation does not
change  the homotopy class of the loop and we can glue the
nullhomotopy to get the required homotopy.
 \QED

\MS

An analogous reasoning may be applied to the case where $M$ is a
nontrivial bundle over $T^3.$ We need first to recall (see section
\ref{sym}) that the space $\mathcal{N}_H \cap L$ has the homotopy
type of the circle. The homotopy described in \ref{mwk} can be
adjusted to keep the base point fixed. Moreover, any closed
1--form whose cohomology class belong to $\mathcal{N}_H \cap L$
can be approximated by a form representing a fibrated class in
$\mathcal{N}_H \cap L$ and which is in addition fibrated. This is
possible since the Chern class $c_1$ is the reduction of an
integer class. Thus we obtain the following

\begin{cor} For a nontrivial circle bundle over $T^3$ we have
$\pi_1S_{inv}=\mathbb Z.$
\end{cor}

\noindent In this context arises a natural problem: does the
Laudenbach--Blank theorem \cite{BL} hold in the parametric case,
i.e. do higher homotopy groups of $\mathcal{N}^{\alpha}$ always
vanish? The answer, as we shall prove below, is negative. In fact,
we will show that for some torus bundles over $S^1$ and any
1--form $\alpha$ defining a fibration, the space $\mathcal{N}^\al$
has the homotopy type of the circle $S^1$.

\MS

\noindent Let $N^3$ be any torus bundle over $S^1$ with the
monodromy given by the matrix $$A=\left(\begin{array}{cc} 1 &k \\
0 & 1 \end{array}\right),$$ where $k$ is a nonzero integer. Let
also $\alpha$ be a closed and nondegenerate 1--form defining a
fibration of $N^3$.

\begin{prop}\label{imp} Under the above assumptions $\mathcal{N}^\alpha$ has
the homotopy type of the circle.
\end{prop}

{\bf Proof}. It is well known that $\alpha$ defines a fibration
with monodromy conjugated to $A$, so we can assume that this
fibration is also defined by $A$.

\noindent We will analyze the fibration

\begin{eqnarray}\label{wr} F{\buildrel \kappa \over \longrightarrow}
Diff_0N^3  \ {\buildrel ev_\alpha \over \longrightarrow} \
\mathcal{N}^\alpha.\end{eqnarray}

\noindent In the present situation we get a fibration $F_0
\rightarrow F {\buildrel \pi' \over \longrightarrow} SO(2)$ such
that the fiber of this fibration consists of those $\phi \in
Diff_0N^3$ which preserve all fibers. Monodromy of this fibration
is given by the formula $\sigma \rightarrow A^{-1}\sigma.$ Next,
$F_0$ can described as a set of isotopies in $T^2$ joining a
diffeomorphism $\phi_0$ with $A\phi_0 A^{-1}.$ For example, in
$F_0$ there are constant paths equal to $A^n$. Taking the map to
$Diff \ T^2$ given by $\phi_0$, we get another fibration over sum
of certain components of $Diff(T^2)$ with fiber equal to
$\Omega_0Diff_0T^2.$ the latter is a contractible apace, thus the
fibration is product. We get $F_0 \cong Diff_0T^2 \times \Omega_0
Diff T^2 \times \mathbb Z \cong T^2 \times \mathbb Z$, where the
last summand is given by powers of $A$. To justify the last
statement we have to prove that any $\phi \in F_0$ induces on a
fixed fiber $T^2$ of $N \ra S^1$ a diffeomorphism isotopic to some
integer power of $A$.

\smallskip

Denote the generators of $H_1(T^2)$ by $a, b$ and the generator of
$H_1(N)$ coming from the base by $\gamma$ (we assume $Aa=a$ and
$Ab=b+ka$ in $H_1(T^2)$). Under such choices the map induced by
the inclusion $i_*:H_1(T^2) \rightarrow H_1(N^3)$ may be
represented by the formula

\begin{eqnarray}<a,b> \rightarrow
<a,b,0> \ \ \text{in} \ \ \mathbb Z \oplus \mathbb Z \rightarrow
\mathbb Z_k \oplus \mathbb Z \oplus \mathbb Z =
H_1(N^3)\end{eqnarray}

\noindent where
\begin{eqnarray}H_1(N^3) \cong <a,b,\gamma \mid
ka> \cong \mathbb Z_k \oplus \mathbb Z^2.\end{eqnarray}

\noindent Any diffeomorphism $\phi \in F_0$ acts trivially on
$H_1(N),$ hence $(\phi_0)_*u \equiv u \mod \ker i_*$ for any $u
\in H_1(T^2).$ This is the case if and only if $$(\phi_0)_*=
\left(\begin{array}{cc} 1 &nk \\
0 & 1 \end{array}\right),$$ since  $\ker i_* =\left\{
\left(\begin{array}{cc} nk \\ 0  \end{array}\right) \mid n \in
\mathbb Z \right\}.$ From this we also see that $$F \cong
Diff_0T^2 \times \Omega_0Diff_0T^2 \times \mathbb R \cong T^2.$$

\noindent In order to examine homotopy properties of $\kappa$ we
will use the theorem of Hatcher \cite{Ha}, which can be stated as
follows. Assume that we have a 3--dimensional manifold $X$ which
is a circle bundle over $T^2$. Denote by $Aut_0(X)$ the set of
automorphisms of the bundle which are isotopic to the identity as
diffeomorphisms of $X$. Then the inclusion $Aut_0(X) \subset
Diff_0(X)$ is a homotopy equivalence.

\MS

\noindent Our manifold $N^3$ can be fibered over $T^2$ with fiber
$S^1$ in the obvious manner, because monodromy $A$ has eigenvalue
1. This bundle has the first Chern class equal to $k$.

\MS

\noindent The fibration

\begin{eqnarray}F_A \rightarrow Aut_0(N)
\rightarrow Diff_0T^2,\end{eqnarray}

\noindent has the fiber $F_A$ equal to the space $Map_0(T^2,S^1)$,
i.e. maps from $T^2$ to $S^1$ homotopic to the constant map. It
can be easily seen by noticing that any $\psi \in Aut_0(N)$
preserves $a$ (which is the fiber of the principal fibration), as
well as $b$ and $\gamma$ (which are some lifts of the generators
of the fundamental group of the base $T^2$). Moreover, if $\psi$
induces a nontrivial loop in $Map(T^2,S^1)$, then it would
transform at least one of the loops $b$ or $\gamma$ to a
nontrivial sum with $a$, which is impossible.

\begin{lemma} The evaluation map \begin{eqnarray} ev :
Map_0(X,S^1) \rightarrow S^1 \ : \ f \mapsto f(*)\end{eqnarray} is
a homotopy equivalence for any compact CW -- complex $X.$
\end{lemma}

{\bf Proof}. Under our assumptions the evaluation map with respect
to a base point $*$ is a fibration with the fiber equal to the
space $Map_0((\cdot,*),(S^1,*))$. But the latter is contractible.
To see this build a map
\begin{eqnarray}Map_0((\cdot,*),(S^1,*)) \rightarrow
Map_0((\cdot,*),(\mathbb{R},*)),\end{eqnarray} \noindent which
assigns to a map $\phi \in Map_0((\cdot,*),(S^1,*))$ its lifting.
This assignment is a homeomorphism, so the thesis follows. \QED

\MS

Consider now the diagram:

$$
\begin{CD}\label{fbi}
             F_A                @>{\hbox{ev}}>>   S^1 \\
             @V{}VV              @V{ }VV           \\
             Aut_0(N)           @>{\hbox{ev}}>>    N \\
             @V{}VV              @V{ }VV           \\
             Diff_0T^2          @>{\hbox{ev}}>>   T^2. \\
\end{CD}
$$

Applying the homotopy exact sequence to the vertical fibrations we
get that $ev:Aut_0(N) \rightarrow N$ induces isomorphisms on all
homotopy groups, so it is a homotopy equivalence.

\smallskip

Let us consider now evaluation for $Diff_0N.$ This gives us the
following commuting diagram

$$
\begin{CD}\label{er}
             F_1           @>{pr}>>  Diff_0T^2      @>{ev}>>     T^2         \\
             @VVV                                                            \\
             F                           @.   @.                   @VVV      \\
             @V{}VV                                                          \\
             Diff_0N      @<<<        Aut_0(N)      @>{ev}>>         N       \\
\end{CD}
$$

where $F_1$ is the component of $id$ in $F_0,$ vertical arrows are
given by inclusion and $pr: F_1 \ra Diff_0T^2$ is given by
restriction to the fixed fiber of the fibration $T^2 \ra N \ra
S^1.$
\smallskip

Since all evaluation maps, as well as $pr$ and the inclusion $F_1
\subset F$ are homotopy equivalences, we see that $F \subset
Diff_0N$ is up to homotopy equal to the inclusion $T^2 \subset N.$
The latter induces monomorphism in $\pi_1$ with coimage isomorphic
to $\mathbb Z.$

\smallskip

From the exact sequence of fibration \ref{wr}

\begin{eqnarray}0 \rightarrow \pi_2\mathcal{N}^\alpha
\rightarrow \pi_1F \rightarrow \pi_1Diff_0N^3 \rightarrow
\pi_1\mathcal{N}^\alpha \rightarrow 0.\end{eqnarray}

\noindent we see now that $\pi_2\mathcal{N}^\alpha$ and higher
homotopy groups of $\mathcal{N}^\alpha$ vanish and
$\pi_1\mathcal{N}^\alpha = \mathbb Z$.  \QED

\smallskip

Let $\mathcal{S}^\alpha \subset S_{inv}$ be the subset of all
invariants symplectic forms whose cohomology class $\alpha \in
H^1(M; \mathbb R)$ is fixed. Imitating the proof of the Lemma
\ref{im} we easily get

\begin{cor} If $M$ is a total space of any circle bundle over $N$
such that $\alpha \in L$ and $\al$ is fibrated, then
$\mathcal{S}^\alpha$ has the homotopy type of the circle.
\end{cor}

\MS \noindent Higher homotopy groups of $S_{inv}$ might be
nontrivial. We shall give here simple example.

\begin{ex}\end{ex} For three--dimensional torus $T^3$ let us define the
map $\Phi:S^2 \rightarrow \mathcal{N}$ by the formula
\begin{eqnarray}\Phi(\theta) = \theta_1 dx_1 + \theta_2 dx_2 +
\theta_3 dx_3,\end{eqnarray} \noindent where $\theta_i$ are
cartesian coordinates of $\theta \in S^2$. For the trivial circle
bundle over $T^3$ the lift of the map $\Phi$ to $S_{inv}$ gives a
nontrivial element of $\pi_2(S_{inv})$.

\MS

Classification of symplectic forms in $S_{inv}$ in case when
$c_1(\xi)$ is finite in $H^2(N^3,\mathbb Z)$ is more simple
because there is a closed connection form on $M^4$.

\section{Almost complex structures}

\noindent We show now how to apply Thurston's theory \cite{Th} to
prove that in the invariant case existence of a path between two
almost complex structures yields a homotopy of symplectic forms
which are tamed by these structures.

\begin{thm}\label{vf} Let $M$ be a closed 4--manifold equipped
with a free circle action. Assume there are two invariant almost
complex structures $J_0$ and $J_1$ on $M^4$ tamed by invariant
symplectic forms $\omega_0$ and $\omega_1$ respectively. If there
exists a path of invariant almost complex structures $\{ J_t \}_{t
\in [0,1]}$ joining $J_0$ and $J_1,$ then $\omega_0$ and
$\omega_1$ belong to the same component of $S_{inv}$ (they are
homotopic through invariant symplectic forms).
\end{thm}

{\bf Proof}. We will use \ref{im} and the Thurston norm. Define
first a family of vector fields $Y_t = \pi_*(J_tX)$ for $t \in
[0,1]$, so $\alpha_i(Y_i)>0$ for $i= 0,1$ for all $x \in M \slash
S^1.$ Our problem reduces then to the following one. Assume we are
given two nondegenerate and closed 1--forms $\alpha_0$ and
$\alpha_1$ on three--dimensional closed manifold together with a
path of nondegenerate (= nowhere vanishing) vector fields $Y_t$
such that $\alpha_i(Y_i)>0$ holds for $i \in \{ 0,1 \}$. The
question is whether $\alpha_0$ and $\alpha_1$ are homotopic
through closed and nondegenerate 1--forms. We claim that the
answer to this question is positive.

\SK

Denote by $e(\alpha) \in H_1(N^3,Z)$  the Euler class, defined as
the Euler class of two--dimensional oriented subbundle
$ker(\alpha)$. Let us recall a theorem \cite{Th} in slightly more
convenient form from \cite{MT}:

\begin{thm}\label{thu} Suppose $\alpha \in H^1(N^3,Z)$ can  be represented by
a fibration over the circle $N^3 \rightarrow S^1$. If $F \subset
H^1(N^3,\mathbb R)$ is the top--dimensional face of the norm ball
$B_T$ containing $\alpha$, then $\phi(e(\alpha))=-1$ for all $\phi
\in F.$
\end{thm}

This theorem has the following corollary.

\begin{cor}\label{pre} Let $\alpha_0$ and $\alpha_1$ be two nondegenerate
and closed 1--forms on closed three--dimensional manifold $N^3$.
If $e(\alpha_0)=e(\alpha_1),$ then $\alpha_0$ and $\alpha_1$ are
homotopic through closed and nondegenerate 1--forms.
\end{cor}

{\bf Proof}. Assume first that the Thurston norm on $N$ does not
vanish identically and that $\alpha_0$ and $\alpha_1$ are in
different faces $F_0$ and $F_1$ respectively. Then, for any $\phi
\in F_0, \psi \in F_1$ we have $\phi(e)=\psi(e)=-1$, where $e$
denotes the class $e(\alpha_0)=e(\alpha_1)$. Thus the functional
$$\eta \mapsto \eta(e)$$ on $H^1(N,\mathbb R)$ takes the value $-1$ on
two different codimension $1$ planes, which is not possible.

If the Thurston norm on $N$ vanishes identically, then it follows
from the assumptions that all nonzero classes possess a closed and
nondegenerate representative, so $\alpha_0$ and $\alpha_1$ are
homotopic by \ref{lab}. \QED

\smallskip

With our assumptions $e(\alpha_0)=e(\alpha_1),$ hence proposition
\ref{pre} implies Theorem \ref{vf}. \QED

It is not true in general that a given path $\{ J_t \}_{t \in
[0,1]}$ lifts, nevertheless there exists another invariant path
$\{ J'_t \}_{t \in [0,1]}$ joining $J_0$ with $J_1$ which can be
lifted. We give now such example.

\begin{ex}\label{to}\end{ex}

It is enough to find a three--dimensional manifold $N^3$ with two
closed, nondegenerate 1--forms $\{\alpha_i\}_{i \in \{0,1\}}$ and
a path $\{ Y_t \}_{t \in [0,1]}$ of nondegenerate vector fields
such that the following conditions are met:

\begin{enumerate}
\item $\alpha_i(Y_i)>0$ for $i \in \{0,1\},$ \item there is no
path $\{ \alpha_t \}_{t \in [0,1]}$ of closed, nondegenerate
1--forms joining $\alpha_0$ and $\alpha_1$ such that
$\alpha_t(Y_t)>0$ for all $t \in [0,1]$ is not satisfied.
\end{enumerate}

Having such data we define $\omega_0,\omega_1,J_t$ on $N^3 \times
S^1$ with desired properties as follows. First, let $B_t$ denote
the subbundle of $TN$ orthogonal (with respect to any Riemannian
metric) to $Y_t.$ It is oriented and of dimension 2, hence there
is a continuous path $J_t^\perp$ of complex structures on $B_t.$
Since the space of almost complex structures (with a given
orientation) on $\mathbb{R}^2$ is contractible we can require the
initial and terminal structures $J_0^\perp$ $J_1^\perp$ to be
prescribed. We extend $J_t^\perp$ to an almost complex structure
on $N^3 \times S^1$ by $JX=Y_t.$ The forms $\alpha_0, \alpha_1$
and $\eta=dx_0$ give invariant symplectic forms $\omega_0,
\omega_1$ as in \ref{sim} tamed by $J_0$ and $J_1$ respectively.

Take for $N$ the three--dimensional torus $T^3$ with the standard
Riemannian metric and define $$\alpha_0=dx_1,\ \alpha_1=dx_2,\
Y_0=\frac{\partial}{\partial x_1}, \ Y_1=\frac{\partial}{\partial
x_2},$$ $$ \omega_0=dx_0 \wedge dx_1 +dx_2 \wedge dx_3,\
\omega_1=dx_0 \wedge dx_2 +dx_3 \wedge dx_1,$$
$$J_0^\perp\left(\frac{\partial}{\partial x_2}\right)=\frac{\partial}{\partial
x_3},\ J_1^\perp\left(\frac{\partial}{\partial
x_3}\right)=\frac{\partial}{\partial x_1},$$ where $x_0$ denotes
the coordinate system on $S^1 \subset S^1 \times N.$ Using the
canonical trivialization of the tangent bundle $TT^3$ we identify,
after normalization, nondegenerate vector fields with maps
$\{T^3\rightarrow S^2\}$. Choose an unknotted circle $$S^1 \subset
D^3\subset T^3.$$ It suffices now to build a map $f:T^3\rightarrow
S^2$ such that $f$ restricted to $S^1$ is equal to a nonzero
vector tangent field to the circle and $f$ is homotopic with the
constant map. The tangent field to the circle can be clearly
extended to a null--homotopy map, so we need only to show that
there does not exist a closed and nondegenerate 1--form $\alpha$
strictly positive on the tangent vector filed to the circle. The
proof is immediate: we have $\int_{S^1}\alpha=0$ due to
contractibility of the circle, while positiveness (condition 1)
requires $\int_{S^1}\alpha>0$. $\Box{}$

\BS

\noindent Theorem \ref{vf} affirmatively answers one of the basic
questions concerning existence and uniqueness problem for
symplectic structures. Namely, the set of invariant symplectic
forms in a fixed homotopy class of invariant almost complex
structures is connected.

\MS

\noindent Consider the following simpleminded question about
existence of symplectic structures: if $J$ is an almost complex
structure, does $M$ admit a symplectic structure in the homotopy
class of $J$? In general the answer is negative, for instance on
$T^4$ there are almost complex structures having nonzero the first
Chern class, while it is zero for those which are compatible with
a symplectic structure. In the invariant case one can give an
answer in terms of the quotient manifold $N$. The obvious
necessary condition is the existence of a closed 1--form $\alpha
\in \mathcal{N}_L$ such that $\alpha(\pi_*JX)>0$ on $N^3 \cong M^4
\slash S^1$. We shall prove that this is also sufficient. Together
with Theorem \ref{vf} it gives in fact a description of
$\pi_0S_{inv}$ as $\pi_0$ of a subspace of invariant almost
complex structures on $M.$

\begin{thm}\label{co} Let $M^4$ be a closed manifold equipped with
a free circle action with infinitesimal generator $X$ and an
invariant almost complex structure $J$. Assume that there exists a
closed and nondegenerate 1--form $\alpha \in \mathcal{N}_L$ such
that $\alpha(\pi_*JX)>0$. Then $J$ is homotopic (through invariant
almost complex structures) to an almost complex structure tamed by
some invariant symplectic form.
\end{thm}

{\bf Proof}. The almost complex structure $J$ defines a complex
structure on the bundle $TM.$ Consider invariant hermitian metric
$<\cdot, \cdot>$ for $J$. Let $A$ be the orthogonal complement of
$X$ on the complex bundle $(TM,J).$ Observe that $A$ is invariant,
$J$--invariant and together with $Lin_{\mathbb R}\{ JX \}$ defines
a connection $\eta$ for the principal circle bundle $M \ra M
\slash S^1 \cong N.$ It enables us to lift $\ker \alpha$ to a
subbundle of $TM$, denoted by $S$. If $S$ coincides with $A,$ we
are essentially done. If not, we have to deform $J$ to obtain this
property. We homotype $J$ to another almost complex structure $J'$
such that

\begin{eqnarray}J'(S)=S \ \text{and} \ J' \upharpoonright
Lin\{X,JX\} = J\upharpoonright Lin\{X,JX\}.\end{eqnarray}

\noindent Since in $\ker\eta$ we have $A \perp \{JX\}$ and (from
the assumption) $S \perp \{JX\}$, we can homotype $J$ to $J'$ by
use the map $S \mapsto A$ given by the projection parallel to
$JX.$ We can proceed for example as follows. The tangent bundle
$TM$ can be decomposed as $Lin_{\mathbb C}\{X\} \oplus S$ and
$Lin_{\mathbb C}\{X\} \oplus A$, so and $J$ orients $S.$
Furthermore whenever we choose unique orthonormal vectors
$U_A=A^\perp$ and $U_S=S^\perp$ in $\ker\eta$ appointed by choices
of orientations a condition $U_A+U_S=0$ can not be satisfied. Then
the shortest path $\{U_t\}_{t \in [0,1]}$ of orthonormal vectors
joining $U_A$ and $U_S$ gives the path $U_t^\perp$ of oriented,
two--dimensional subbundles in $\ker\eta$. We induce almost
complex structures on each $U_t^\perp$ from $A$ using the parallel
projection $A \mapsto U_t$ to $JX.$

\SK

\NI Let $\beta$ be a closed 2--form on $N$ such that $\beta \wedge
\alpha$ is nowhere vanishing , i.e. $\beta \upharpoonright
\ker\alpha$ nowhere vanishes. We immediately see that
$\beta(\pi_*W,\pi_*J'W)$ is of constant sign for nonzero $W \in
S$, so replacing $\beta$ by $-\beta$ if necessary we have assume
that $\beta(\pi_*W,\pi_*J'W)>0.$

\noindent We now use \ref{sim} to define an invariant form by
\begin{eqnarray}\omega = \eta \wedge \pi^*\alpha + \pi^*(K
\beta+\phi),\end{eqnarray} where $K$ big enough to guarantee that
$\omega$ is symplectic and  satisfy the condition
\begin{eqnarray}\pi^*(K \beta+\phi)(W,J'W)>0\end{eqnarray} for
all nonzero $W \in S$.

\SK

\NI The form $K\beta+\phi$ has maximal rank, so the two--sided
kernel $$\{u \mid \forall v \ \ (K\beta+\phi)(u,v)=0\}$$ is a one
dimensional subbundle $T$ of $TN$. Furthermore the bundle $T$ is
oriented because it is transversal to the canonically oriented
bundle $\ker \alpha$ due to condition
$\alpha(\pi_*J'X)=\alpha(\pi_*JX)>0$. We are able to choose a
section $Z$ of the bundle $T$ such that $\alpha(Z)>0$. We want now
to homotype $J'$ to $J''$ in such a manner so the conditions
\begin{eqnarray}J''(S)=S \ \text{and} \ \pi_*J''X=Z\end{eqnarray}
hold. The vector fields $Z$ and $\pi_*J'X$ are both positive on
$\alpha$ and so is the convex path $\{Z_t\}_{t \in [0,1]}$ joining
these fields. By transversality $\ker\alpha$ and $Z$ we get that
$\{\widetilde{Z_t}\}_{t \in [0,1]}$ are transversal to $S$ for all
$t \in [0,1]$. Naturally, we define $$J'_tX=\widetilde{Z_t}$$
where $J'_0=J'$ and $J'_1=J''.$

\SK

\noindent We will show that $J''$ is an invariant almost complex
structure tamed by $\omega$.

\smallskip

\noindent Observe that $$\omega(X,J''X)=\eta \wedge
\pi^*\alpha(X,J''X)>0$$ and $$\omega(V,J''V)=\pi^*(K \beta+\phi)
(V,J''V)>0$$ for any nonzero $V \in S$. Take now any nonzero $W
\in TM$. Using the decomposition $T_xM =Lin_{\mathbb R}\{X,J''X\}
\oplus S$ we can check directly that $J''$ is tamed by $\omega$.
\QED

We give now an example of invariant almost complex structure on
$T^4$ which has no taming invariant symplectic form, hence does
not satisfy the assumption of \ref{co}. Almost complex structures
on $T^4$ correspond to maps $T^4 \ra GL(4,\mathbb R) \slash
GL(2,\mathbb C) \cong S^2.$ Take the standard trivialization of
the tangent bundle $TT^4$ and the standard (i.e. corresponding to
the constant map) almost complex structure. Twist $J$ on $D^3
\times S^1 \subset T^3 \times S^1$ to correspond to the map
$$T^4 \ra D^3 \times S^1 \slash \partial(D^3 \times S^1) \ra D^3 \slash
\partial D^3 {\buildrel \hbox{Hopf} \over
\longrightarrow} S^2.$$

There is no symplectic form which tames the resulting $J.$ It
follows from \ref{tor} since any two invariant symplectic forms on
$T^4$ are homotopic. By the same argument any invariant almost
complex structure nonhomotopic to the standard one cannot be tamed
by an invariant symplectic form.

\SK

\newcommand{\gsv}{\mbox{$\tilde{\GG}_s(V)$}}
\newcommand{\gsw}{\mbox{$\tilde{\GG}_s(W)$}}
\newcommand{\gso}{\mbox{$\tilde{\GG}_s(0)$}}
\newcommand{\gsvw}{\mbox{$\tilde{\GG}_s(V\times W)$}}

\bibliographystyle{amsalpha}

\medskip
\noindent {\bf Mathematical Institute, Wroc\l aw University,

\noindent pl. Grunwaldzki 2/4,

\noindent 50--384 Wroc\l aw, Poland}

\begin{flushleft}
\tt hajduk@math.uni.wroc.pl

\tt rwalc@math.uni.wroc.pl
\end{flushleft}

\end{document}